\documentclass[11pt]{article}

\usepackage{latexsym,amssymb}
\usepackage{enumerate}
\usepackage{amsmath}
\usepackage[latin1]{inputenc}
\usepackage[english]{babel}

\textheight20cm
\textwidth14cm
\oddsidemargin1cm


\newcommand{\bg}{\begin}
\newcommand{\e}{\end}
\newcommand{\be}{\bg{enumerate}}
\newcommand{\ee}{\e{enumerate}}
\newcommand{\bi}{\bg{itemize}}
\newcommand{\ei}{\e{itemize}}
\newcommand{\ba}{\begin{array}}
\newcommand{\ea}{\end{array}}

\newcommand{\no}{\noindent}
\newcommand{\sub}{\subseteq}

\newcommand{\Ra}{\Rightarrow}

\newcommand{\lgra}{\longrightarrow}


\newcommand{\MV}[1]{\langle #1,\oplus,\lnot,0 \rangle}

\newcommand{\bfm}[1]{\mbox{\boldmath$#1$}}


\title{ Locally finite quasivarieties of MV-algebras.\thanks{This work is
partially supported by Grants  SGR/2000--0007 of D.G.R. of
Ge\-neralitat de Catalunya and by Grant PB97-0888 of D.G.I.C.Y.T. of Spain.}}

\author{ {\sc Joan Gispert \  and \ Antoni Torrens}}
\date{}

\bg{document}
\maketitle

\bg{abstract}
In this paper we show that every locally finite quasivariety of MV-algebras
is finitely generated and finitely based. To see this result we study
critical MV-algebras. We also give  axiomatizations of some of these
quasivarieties.

\e{abstract}

\section*{Introduction}

In \cite{Ch1,Ch2}  C.C.Chang introduced MV-algebras in order to
give an algebraic counterpart of the  \L ukasiewicz's many valued
propositional calculus. In fact, the class of all MV-algebras, in a termwise
equivalent presentation named Wajsberg algebras, is the equivalent variety
semantics, in the sense of \cite{B-P}, of this calculus (see \cite{RTV}).

From the equivalence between the class of MV-algebras and
\L ukasiewicz logic,  it is easy to see that finitary extensions
of \L ukasiewicz's propositional calculus  correspond to subquasivarieties of
MV-algebras, and axioms and rules of the calculus  correspond with equations and
quasiequations, respectively. Hence, finite axiomatizable finitary  extensions
of \L ukasiewicz's propositional calculus correspond with finite axiomatizable
quasivarieties of MV-algebras.

In this paper, we study  finite axiomatizability of locally finite quasivarieties
of MV-algebras. To be precise, we  show in Section 2 that locally finite
quasivarieties and finitely generated quasivarieties of MV-algebras coincide
(Theorem \ref{locfinfingen}) and that they are finitely axiomatizable
(Theorem
\ref{finitaxiomatitzable}). To prove these results,  we give a  characterization
of  critical MV-algebras (Theorem \ref{Wcritical}) and  we see that any
locally finite quasivariety of   MV-algebras is generated by  critical
MV-algebras (Theorem \ref{Criticalgen}).

Finally, in Section 3, we give two examples of  locally finite quasivarieties with their
finite axiomatizations.

We include a preliminary section, Section 1, containing   basic
definitions, results and notation used in the paper.

\section{ Definitions and first properties.}

 An \textbf{ MV-algebra} is an algebra $\mathbf{ A}=\langle A,\oplus,\lnot,0
\rangle$  of type $(2,1,0)$ satisfying the following equations:
\bg{description}
 \item\textbf{ MV1.} $(x\oplus y)\oplus z\approx x\oplus
  (y\oplus z)$
 \item\textbf{ MV2.} $x\oplus y\approx
  y\oplus x$
 \item\textbf{ MV3.} $x\oplus 0\approx x$
 \item\textbf{ MV4.} $\lnot(\lnot x)\approx x$
 \item\textbf{ MV5.} $x\oplus \lnot 0\approx \lnot 0$
 \item\textbf{ MV6.} $\lnot (\lnot
  x\oplus  y)\oplus y\approx \lnot(x\oplus \lnot y)\oplus
  x$
\end{description}
By taking $y=\lnot 0$ in MV6, we deduce:
\bg{description}
 \item\textbf{ MV7.} $x\oplus \lnot x\approx \lnot 0$.
\end{description}
More precisely, if we set $1=\lnot 0$ and
$x\odot y=\lnot(\lnot x\oplus\lnot y)$, then
$\langle A,\oplus,\odot,\lnot,0,1 \rangle$ satisfies all
axioms given in \cite[Lemma 2.6]{Mu1}, and hence the above
definition of MV-algebras is equivalent to Chang's
definition \cite{Ch1}. We denote the class of all MV-algebras by $\mathbb{W}$.
Since it is an equational class, $\mathbb{W}$ is a variety.

Given  a {\em lattice ordered abelian group} $\mathbf{ G}=\langle
G,\land ,\lor ,+,-,0\rangle$   and $u\in G\quad u>0$, we define the algebra $\Gamma(\mathbf{
G},u)=\MV{[0, u]}$ where
\bi
 \item[-] $[0, u]=\{a\in G : 0\leq a\leq u\}$,
 \item[-] $a\oplus b=u\land(a+b)$, $ \neg a=u-a$ and $ 0=0^\mathbf{ G}$.
\ei
Then $\Gamma(\mathbf{ G},u)$ is an MV-algebra. In fact any MV-algebra is isomorphic
to the unit segment of some lattice ordered abelian group. Concretely, the
category of MV-algebras is equivalent to the category of lattice ordered
abelian groups with strong unit (See  \cite{CDM},\cite{Mu1}).

The following MV-algebras play an important role in the paper.
\bg{itemize}
\item $\displaystyle\mathbf{ [0,1]}=\Gamma(\mathbf{ R},1),$ where $\mathbf{ R}$ is the totally
ordered group of the reals.
\item $\displaystyle \mathbf{ [0,1]}\cap\mathbf{ Q}=\Gamma(\mathbf{
  Q},1)=\MV{\{\frac{k}{m} : k\leq m<\omega\}}$,
  where $\mathbf{ Q}$ is the totally ordered abelian group of the
  rationals and $\omega $ represents the set of all natural numbers.
\ei
For every $0<n<\omega $
\bi
\item $\displaystyle \mathbf{ \L}_{n}=\Gamma(\mathbf{ Q}_{n},1)=\MV{\{\frac{k}{n} :
 0\leq k\leq n\}}$,
 where $\displaystyle \mathbf{ Q}_{n}=\{\frac{k}{n} : k\in\mathbb{Z}\}$ is a
 subgroup of $\mathbf{Q}$ and $\mathbb{Z}$ is the set of all integers.
\item  $\displaystyle \mathbf{ \L}_{n}^{\omega}=\Gamma(\mathbf{ Q}_{n}\otimes
 \mathbf{ Z},(1,0))=\MV{\{(\frac{k}{n},i) :(0,0)\leq(\frac{k}{n},i)\leq (1,0)\}}$,
 where $\mathbf{Z}$ is the totally ordered group of the integers and
 $\mathbf{ Q}_{n}\otimes\mathbf{ Z}$  is the lexicographic product of
 $\mathbf{ Q}_{n}$ and  $\mathbf{ Z}$.
\e{itemize}

The following theorem states some well-known results on   simple and/or  finite
MV-algebras.( See for instance \cite{CDM}).
\bg{theorem}
\label{Coneguts} \mbox{ }
\be
 \item Every simple MV-algebra is isomorphic to a subalgebra of \textbf{ [0,1]}.
 \item Every finite simple MV-algebra is isomorphic to $\mathbf{ \L}_{n}$ for some
  $n\in\omega $.
 \item Every finite MV-algebra is isomorphic to a direct product
  of finite simple MV-algebras.
 \item $\mathbf{ \L}_{n}\sub\mathbf{ \L}_{m}$ if and only if $n|m$. \hfill $\Box$
\ee
\e{theorem}

We denote by  $\mathbb{I}$, $\mathbb{H}$, $\mathbb{S}$, $\mathbb{P}$,
$\mathbb{P}_{R}$ and $\mathbb{P}_{U}$ the   operators {\em isomorphic
image, homomorphic image, substructure, direct product,
reduced product and ultraproduct} res\-pectively. We recall
that a class $\mathbb{K}$ of algebras is a \textbf{ variety} if and only if it is
closed under $\mathbb{H}$, $\mathbb{S}$ and $\mathbb{P}$. And a class
$\mathbb{K}$ of algebras is a \textbf{ quasivariety} if and only if it is
closed under $\mathbb{I}$, $\mathbb{S}$ and $\mathbb{P}_{R}$, or equivalently, under
$\mathbb{I}$, $\mathbb{S}$, $\mathbb{P}$ and $\mathbb{P}_{U}$.
Given a class $\mathbb{K}$ of algebras, the variety generated by ${\Bbb
K}$, denoted by $\mathbb{V}(\mathbb{K})$, is the least variety containing
$\mathbb{K}$. Similarly, the quasivariety generated by a class
$\mathbb{K}$, which we denote by $\mathbb{Q}(\mathbb{K})$, is the least quasivariety
containing $\mathbb{K}$. We also recall that a class $\mathbb{K}$ of algebras is a
variety if and only if it is an equational class, and $\mathbb{K}$ is a
quasivariety if and only if it is a quasiequational class.

The subvarieties of $\mathbb{W}$ are well known:
\bg{theorem}{\rm \cite[Theorem 4.11]{Ko}}
\label{varietats} $\mathbb{K}$ is a proper subvariety of $\mathbb{W}$ if and only if
there exist two disjoint finite subsets $I,J$ of natural numbers such that
\[\mathbb{K}=\mathbb{V}(\mathbf{ \L}_{i}\ i\in I,\mathbf{ \L}_{j}^{\omega}\ j\in J).\qquad \qquad \Box\]
\e{theorem}

An algebra $\mathbf{ A}$ is \textbf{ locally finite} if and only if every finitely
generated subalgebra is finite. A class $\mathbb{K}$ is \textbf{ locally finite} if
and only if every member of $\mathbb{K}$ is locally finite.\\
 If $\mathbb{V}$ is a variety, then
$\mathbb{V}$ is locally finite if and only if all free algebras with respect $\mathbb{V}$
over a finite set of generators are finite \cite[page 69]{Bu}.\\
If $\mathbb{K}$ is a quasivariety, then $\mathbb{K}$ is locally finite if and only if ${\Bbb
K}$ is contained in a locally finite variety.

A variety, or a quasivariety, is \textbf{ finitely generated} if it is
generated by a finite set of finite algebras.

\bg{theorem}{\em \cite[page 70]{Bu}}
\label{varlocfin} Every finitely generated variety is locally finite. \hfill $\Box$
\e{theorem}

From the above we have:
\bg{lemma}
\label{varWlocfin} Let $\mathbb{K}$ be a variety of MV-algebras. $\mathbb{K}$ is
a locally finite variety if and only if $\mathbb{K}=\mathbb{V}(\mathbf{
\L}_{n_{1}},\ldots\mathbf{ \L}_{n_{r}})$ for some $n_{1},\ldots,n_{r}\in\omega$
\e{lemma}
\textbf{ Proof :} By Theorem \ref{varlocfin}, for every  $r<\omega$ and any
$n_{1},\ldots,n_{r}\in\omega$,  $\mathbb{V}(\mathbf{ \L}_{n_{1}},\ldots\mathbf{
\L}_{n_{r}})$  is a locally finite variety.


\noindent Since  $\mathbf{ \L}_{1}^{\omega}\in\mathbb{W}$ is infinite and it is finitely generated
by the element $(0,1)$, $\mathbb{W}$ is not locally finite.\\
Assume $\mathbb{K}$ is a proper variety of $\mathbb{W}$. If  $\mathbb{K}$ is not
of the form $\mathbb{V}(\mathbf{ \L}_{n_{1}},\ldots\mathbf{ \L}_{n_{r}})$ for some
$n_{1},\ldots,n_{r}\in\omega$, then by Theorem \ref{varietats} we have that
$\mathbb{K}=\mathbb{V}(\mathbf{ \L}_{i}\; i\in I,\mathbf{
\L}_{j}^{\omega}\; j\in J)$ with $J\neq\emptyset$. Hence there is $j\in J$
such that $\mathbf{ \L}_{j}^{\omega}\in \mathbb{K}$. Since $\mathbf{ \L}_{j}^{\omega}$ is infinite
and it is finitely generated by $\displaystyle (0,1),
(\frac{1}{j},0)\in {\rm \L}_{j}^{\omega}$, $\mathbf{ \L}_{j}^{\omega}$ is not locally finite
and therefore, $\mathbb{K}$ is not a locally finite variety.
\hfill
$\Box$

\section{ Locally finite quasivarieties and critical algebras.}

We  want to obtain all locally finite
 quasivarieties of
MV-algebras. First we observe that from Theorem \ref{varlocfin} it follows.
\bg{corollary}\label{FinGenisLoc}
 Every finitely generated quasivariety is locally
finite.\hfill $\Box$
\e{corollary}
And from Lemma \ref{varWlocfin} we obtain:
\bg{corollary}
\label{quasivarWlocfin} A quasivariety of MV-algebras is locally finite if and
only if it is a subquasivariety of a variety of the form $\mathbb{V}(\mathbf{
\L}_{n_{1}},\ldots\mathbf{ \L}_{n_{r}})$ for some
$n_{1},\ldots,n_{r}\in\omega$.\hfill $\Box$
\e{corollary}
A \textbf{ critical} algebra is a finite algebra not belonging to the
quasivariety generated by all its proper subalgebras.
The interest of critical algebras is given by the following
result, which is mentioned  in  \cite[page 128]{Dz0}, but no proof is given.
\bg{theorem}
\label{Criticalgen}  Every locally finite quasivariety is generated by its
critical algebras.
\e{theorem}
\textbf{ Proof :} Let $\mathbb{K}$ be a locally finite quasivariety  and $\{\mathbf{ A}_{i} : i\in I\}\sub
\mathbb{K}$ the family of all critical algebras contained in $\mathbb{K}$. Obviously, ${\Bbb
Q}(\{\mathbf{ A}_{i} : i\in I\})\sub\mathbb{K}$. Assume $\displaystyle
\mathop{\&}\limits_{k=1}^{m}\varphi_{k}\approx\psi_{k}\Ra\varphi\approx\psi(x_{0},\ldots
x_{n})$ is a quasiequation not satisfied by $\mathbb{K}$. Therefore, there
exist an algebra $\mathbf{ A}\in \mathbb{K}$ and $a_{0}\ldots,a_{n}\in A$ such that
\[\mathbf{ A}\not\models \mathop{\&}\limits_{k=1}^{m}
\varphi_{k}\approx\psi_{k}\Ra\varphi\approx\psi(a_{0},\ldots, a_{n})\]
Let ${\bfm [\![a_{0}\ldots,a_{n}]\!]_{\mathbf{ A}}}$ be the subalgebra of $\mathbf{
A}$ generated by $\{a_{0},\ldots a_{n}\}$. Since $\mathbb{K}$ is locally finite,
\bi
 \item[a)] ${\bfm [\![a_{0}\ldots,a_{n}]\!]_{\mathbf{ A}}}$ is finite
 \item[b)] ${\bfm [\![a_{0}\ldots,a_{n}]\!]_{\mathbf{ A}}}\not\models
  \mathop{\&}\limits_{k=1}^{m}
  \varphi_{k}\approx\psi_{k}\Ra\varphi\approx\psi$
\ei
It is easy to see that a quasivariety generated by a finite algebra $\mathbf{B}$
is generated by all  critical subalgebras of \textbf{B} (The proof is
straightforward by induction over the cardinal of \textbf{B}). Therefore,
from a)  and b)  we deduce that  there is a critical algebra
$\mathbf{A}_{l}\sub{\bfm [\![a_{0}\ldots,a_{n}]\!]_{\mathbf{ A}}}$ such that
\[\mathbf{ A}_{l}\not\models \mathop{\&}\limits_{k=1}^{m}
\varphi_{k}\approx\psi_{k}\Ra\varphi\approx\psi\]
So, since $\mathbf{ A}_{l}\in\{\mathbf{ A}_{i} : i\in I\}$ and a quasivariety is a
class of algebras definable by means of quasiequations, $\mathbb{K}\sub\mathbb{Q}(\{\mathbf{
A}_{i} : i\in I\})$ \hfill
$\Box$

Our next purpose is to characterize critical MV-algebras. We need a previous
result.
\bg{lemma}
\label{Ln12inclusiolm12} If $\mathbf{ \L}_{n_{0}}\times\cdots\times\mathbf{
\L}_{n_{l-1}}$ is embeddable into $\displaystyle \prod_{j\in J}\mathbf{ \L}_{m_{j}}$
 where the set \break $\{m_{j} : j\in J\}$ is finite then
\be
 \item For every $i<l$ there exists $j\in J$ such that $n_{i}|m_{j}$.
 \item For every $j\in J$ there exists $i<l$ such that $n_{i}|m_{j}$.
\ee
\e{lemma}
\textbf{Proof :}
1) If $\mathbf{ \L}_{n_{0}}\times\cdots\times\mathbf{ \L}_{n_{l-1}}$ is  embeddable
  into $\displaystyle \prod_{j\in J}\mathbf{ \L}_{m_{j}}$,  then
  \[ \mathbf{ \L}_{n_{1}}\times\cdots\times\mathbf{ \L}_{n_{l}}\in {\Bbb
  V}(\prod_{j\in J}\mathbf{ \L}_{m_{j}}) =\mathbb{V}(\{\mathbf{ \L}_{m_{j}}; j\in J\}).\]
  Hence, for every $i<l$, $\mathbf{ \L}_{n_{i}}\in\mathbb{V}(\{\mathbf{
  \L}_{m_{j}}; j\in J\})$. Since $\{m_{j} : j\in J\}$ is finite,  from a
  result due to J\'onsson \cite[page 149]{Bu}, we deduce that the class subdirectly
  irreducible members of $\mathbb{V}(\{\mathbf{ \L}_{m_{j}}; j\in J\})$ is
$\mathbb{I}(\{\mathbf{
  \L}_{n} : \exists j\in J\ \mathbf{ \L}_{n}\sub\mathbf{ \L}_{m_{j}}\})$. Since $\mathbf{
  \L}_{n_{i}}$ is simple, therefore subdirectly irreducible, for every
  $i<l$ there exists $j\in J$ such that $\mathbf{ \L}_{n_{i}}\sub\mathbf{
  \L}_{m_{j}}$,  and by 4 of Theorem \ref{Coneguts} $n_{i}|m_{j}$.\\
2) For each $j\in J$ consider the natural projection:
  $\displaystyle \pi_{j}:\prod_{j\in J}\mathbf{ \L}_{m_{j}}\lgra\mathbf{ \L}_{m_{j}}$.
  Let  $\gamma\colon\mathbf{
  \L}_{n_{0}}\times\cdots\times\mathbf{ \L}_{n_{l-1}}\to\displaystyle
  \prod_{j\in J}\mathbf{ \L}_{m_{j}}$ be an  embedding, then for every $j\in J$
  $\gamma_{j}=\pi_{j}\circ\gamma $ is an homomorphism from $\mathbf{
  \L}_{n_{0}}\times\cdots\times\mathbf{ \L}_{n_{l-1}}$ to $\mathbf{
  \L}_{m_{j}}$. Hence
  \[\mathbf{ \L}_{n_{0}}\times\cdots\times\mathbf{ \L}_{n_{l-1}}/Ker(\gamma_{j})\cong
  \gamma_{j}(\mathbf{ \L}_{n_{0}}\times\cdots\times\mathbf{ \L}_{n_{l-1}})\sub\mathbf{
  \L}_{m_{j}}\]
  So, $\mathbf{ \L}_{n_{0}}\times\cdots\times\mathbf{ \L}_{n_{l-1}}/Ker(\gamma_{j})$
  is  simple, and  by \cite[Theorem 4.1.19]{CDM} we
  have that $Ker(\gamma_{j})$  is a maximal congruence relation of $\mathbf{
 \L}_{n_{0}}\times\cdots\times\mathbf{ \L}_{n_{l-1}}$.
 From \cite[Lemma 2.3]{C-T}  (see also \cite{Rod})  and the fact that all $\mathbf{
\L}_{n}$'s are simple, it can be deduced that  there is $k<l$ such that
\[Ker(\gamma_{j})=
\L_{n_{0}}^{2}\times\cdots\times\L_{n_{k-1}}^{2}\times
\Delta _{\mathbf{ \L}_{n_{k}}}\times\L_{n_{k+1}}^{2}\times\cdots\times
\L_{n_{l-1}}^{2}.\]
 Hence, for every $j\in J$ there exists $k< l$ such that
\[\mathbf{ \L}_{n_{0}}\times\cdots\times\mathbf{
  \L}_{n_{l-1}}/Ker(\gamma_{j})\cong\mathbf{ \L}_{n_{k}}\sub\mathbf{ \L}_{m_{j}}.\] Thus
$n_{k}|m_{j}$. \hfill
$\Box$\\
Finally we give a characterization of all critical MV-algebras.
\bg{theorem}
\label{Wcritical}  An MV-algebra $\mathbf{ A}$ is critical if and only if $\mathbf{
A}$ is isomorphic to a finite MV-algebra $\mathbf{ \L}_{n_{0}}\times\cdots\times\mathbf{
\L}_{n_{l-1}}$ satisfying the following conditions:
\be
 \item For every $i,j<l$, $ i\neq j$ implies $ n_{i}\neq n_{j}$.
 \item If there exists $n_{j}$ $j<l$ such that $n_{i}|n_{j}$ for some
  $i\neq j$, then $n_{j}$ is unique.
\ee
\e{theorem}
\textbf{Proof :}
Assume that  $\mathbf{ A}=\mathbf{ \L}_{n_{0}}\times\cdots\times\mathbf{ \L}_{n_{l-1}}$
satisfies conditions $1)$ and $2)$. First, we will show the following:\\
\textbf{Claim:} {\em Every proper subalgebra of $\mathbf{ A}$ is embeddable into a
subalgebra of
$\mathbf{ A}$ of the form
$\mathbf{ \L}_{d_{0}}\times\cdots\times\mathbf{ \L}_{d_{l-1}}$, where
$d_{i}|n_{i}$ for each $i<l$ and there exists $j< l$ such
that $d_{j}\neq n_{j}$.}\\
\textbf{Proof of the claim:} Let \textbf{B} be a proper subalgebra of \textbf{A}. Since
\textbf{A} is finite \textbf{B} is also finite and by Theorem \ref{Coneguts},
 \textbf{B} is isomorphic to
$\mathbf{ \L}_{p_{0}}\times\cdots\times\mathbf{ \L}_{p_{r-1}}$. For each $i<l$
consider the natural projection:
$\displaystyle \pi_{i}:\mathbf{ A}\to\mathbf{ \L}_{n_{i}}$,  if for
all $i< l$, we write $\gamma_{i}=\pi_{i}\restriction_\mathbf{ B} $, then
we can assume that  $\mathbf{ B}$ is
embeddable into $\gamma_{0}(\mathbf{ B})\times\cdots\times\gamma_{l-1}(\mathbf{ B})$.
Moreover, since
$\gamma_{i}(\mathbf{ B})\subseteq\mathbf{ \L}_{n_{i}}$,  we have
$\gamma_{i}(\mathbf{ B})=\mathbf{ \L}_{d_{i}}$ for some $d_{i}|n_{i}$.\\
Assume that $\gamma_{i}(\mathbf{ B})=\mathbf{ \L}_{n_{i}}$ for each $i< l$.
Then, for every $i< l$, $\mathbf{ B}/Ker(\gamma_{i})\cong \mathbf{
\L}_{n_{i}}$ and  since $\mathbf{ \L}_{n_{i}}$ is simple, $Ker(\gamma_{i})$ is a
maximal congruence relation of $\mathbf{B}$. From \cite[Lemma 2.3]{C-T}  (see also
\cite{Rod})  and the fact that all $\mathbf{ \L}_{n}$'s are simple, it can be deduced
that  there is $k<r$ such that
\[Ker(\gamma_{i})=
\L_{p_{0}}^{2}\times\cdots\times\L_{p_{k-1}}^{2}\times
\Delta _{\mathbf{ \L}_{p_{k}}}\times\L_{p_{k+1}}^{2}\times\cdots\times
\L_{p_{r-1}}^{2}.\]
 Hence, for every $i< l$ there exists $k< r$ such that
$\mathbf{ \L}_{p_{k}}=\mathbf{ \L}_{n_{i}}$. By condition $(1)$,
$i\not= j$ implies $ n_{i}\neq n_{j}$ , so $l\le r$ and
\[\mathbf{ B}\cong\mathbf{
\L}_{n_{0}}\times\cdots\times\mathbf{ \L}_{n_{l-1}}\times\mathbf{
\L}_{m_{l}}\times\cdots\times\mathbf{ \L}_{m_{r-1}} = \mathbf{ A}\times\mathbf{
\L}_{m_{l}}\times\cdots\times\mathbf{ \L}_{m_{r-1}}\]
that implies $|\mathbf{ A}|\leq
|\textbf{B}|$, which  contradicts that   \textbf{B} is a proper subalgebra of
\textbf{A}. And the claim is proved.\\
 Suppose that
$\mathbf{ A}
\in\mathbb{Q}(\{\mathbf{ B} \subsetneq \mathbf{ A}\})$, then
\[\mathbf{ A} \in\mathbb{I}\mathbb{S}\mathbb{P}\mathbb{P}_{U}
(\{\mathbf{ \L}_{d_{0}}\times\cdots\times\mathbf{ \L}_{d_{l-1}}:\  \forall i\
d_{i}|n_{i};\ \exists  k\  d_{k}\neq n_{k}\}).\]
Since $\{\mathbf{ \L}_{d_{0}}\times\cdots\times\mathbf{ \L}_{d_{l-1}}:\;  \forall i\;
d_{i}|n_{i};\; \exists  k\;  d_{k}\neq n_{k}\}$ is a finite set of finite
MV-algebras, we have that $ \mathbf{ A}\in\mathbb{I}\mathbb{S}\mathbb{P}(\{\mathbf{
\L}_{d_{0}}\times\cdots\times\mathbf{
\L}_{d_{l-1}}:\;  \forall i\; d_{i}|n_{i};\; \exists  k\;  d_{k}\neq n_{k}\})$.
Thus  $\mathbf{ \L}_{n_{0}}\times\cdots\times\mathbf{ \L}_{n_{l-1}}$ is embeddable into
$\displaystyle \prod_{k< n}(\mathbf{ \L}_{d_{0,k}}\times\cdots\times\mathbf{
\L}_{d_{l-1,k}})^{\alpha_{k}}$ where
 \[\{\mathbf{ \L}_{d_{0,k}}\times\cdots\times\mathbf{
\L}_{d_{l-1,k}}: k< n\}\sub\{\mathbf{ \L}_{d_{0}}\times\cdots\times\mathbf{
\L}_{d_{l-1}}:\;  \forall i\; d_{i}|n_{i};\; \exists  k\;  d_{k}\neq n_{k}\}.\]
Since the set $\{d_{t,k}:  t< l; k\le n\}$ is finite, we can apply
Lemma \ref{Ln12inclusiolm12}.\\
If there  exists $ i,j<l$ such that $i\neq j$ and $n_{i}|n_{j}$, then  by
condition $2)$, $n_{j}$ is unique.  By Lemma \ref{Ln12inclusiolm12},
there exists $\mathbf{ \L}_{d_{t,m}}$ such that $n_{j}|d_{t,m }$. That is,
 there exists
\[\mathbf{ \L}_{d_{0,m}}\times\cdots\times\mathbf{ \L}_{d_{l-1,m}}
\in\{\mathbf{ \L}_{d_{0}}\times\cdots\times\mathbf{ \L}_{d_{l-1}}:\;  \forall i\;
d_{i}|n_{i};\; \exists  k\;  d_{k}\neq n_{k}\}\]
such that $n_{j}|d_{t,m}$ for some $ t<l$.
By condition $2)$, $n_{j}$ does not divide any $n_{i}$ other than itself.
Therefore, since $d_{t,m}$ is a divisor of $n_{t}$, we have that
$n_{t}=d_{t,m}=n_{j}$ and by $1)$, $t=j$. Since
\[\mathbf{ \L}_{d_{0,m}}\times\cdots\times\mathbf{ \L}_{d_{l-1,m}}\in
\{\mathbf{ \L}_{d_{0}}\times\cdots\times\mathbf{ \L}_{d_{l-1}}:\;  \forall i\;
d_{i}|n_{i};\; \exists  k\;  d_{k}\neq n_{k}\},\]
there exists $ r\not= j< l$ such that $d_{r,m}|n_{r}$ and $d_{r,m}\neq n_{r}$.
By $2)$ of  Lemma \ref{Ln12inclusiolm12}, there exists $s,r<l$ such that $
s\neq r < l$ and $n_{s}|d_{r,m}$. Thus  $n_{s}|n_{r}$, $n_{i}|n_{j}$, $r\neq
s$, $i\neq j$  and $r\neq j$, which contradicts condition $2)$.\\
If for all $1\leq i,j\leq l$ such that $i\neq j$, $n_{i}\not\! | n_{j}$, then
the same argument follows by taking  any $n_{j}$, $j< l$.\\
Since $\mathbf{ A}$ is finite and  $\mathbf{ A}\not\in \mathbb{Q}(\{\mathbf{ B}\subsetneq
\mathbf{ A}\})$, $\mathbf{ A}$ is critical.\\
 Conversely,  if  $\mathbf{ A}$ is a critical MV-algebra, then  $\mathbf{A}$ is
finite and  by Theorem \ref{Coneguts}, we can suppose, without loss of
generality, that
\[\mathbf{ A}={\mathbf{ \L}_{n_{0}}}^{m_{0}}\times\cdots\times{\mathbf{
\L}_{n_{k-1}}}^{m_{k-1}}\] for some $n_{0},\ldots,n_{k-1},m_{0},\ldots,m_{k-1}
\in\omega$ and $n_{i}\neq n_{j}$ when $i\neq j$. If not all $m_{i}$'s are equal
to
$1$, then  the correspondence
\[\alpha \colon(a(0),\ldots,a(k-1))\mapsto
\alpha(a)=\Big((\overbrace{a(0),\ldots,a(0)}^{m_{0}}),\ldots,
(\overbrace{a(k-1),\ldots,a(k-1)}^{m_{k-1}})\Big)\]
defines an isomorphism from
$\mathbf{ \L}_{n_{0}}\times\cdots\times\mathbf{ \L}_{n_{k-1}}$
onto a proper subalgebra of  $\mathbf{ A}$.
Let $m=max\{m_{0},\ldots,m_{k-1}\}$, then the  correspondence
\[\beta \colon{ \L_{n_{0}}}^{m_{0}}\times\cdots\times{
\L_{n_{k-1}}}^{m_{k-1}}\to \L_{n_{0}}^{m}\times\cdots\times
\L_{n_{k-1}} ^{m}\]
such that for every  $r< k$,
\[\beta ((b(0),\ldots,b(k-1)))(r)=(b(r)(1),\ldots,b(r)(m_{r}),
\overbrace{b(r)(1), \ldots,b(r)(1)}^{m-m_{r}}) ,\]
gives an embedding from $\mathbf{ A}$ into
${\mathbf{ \L}_{n_{0}}}^{m}\times\cdots\times{\mathbf{
\L}_{n_{k-1}}}^{m}\cong \big( \mathbf{ \L}_{n_{0}}\times\cdots\times\mathbf{
\L}_{n_{k-1}} \big)^{m}.$
Thus $\mathbf{ A}\in \mathbb{Q}(\mathbf{ \L}_{n_{0}}\times\cdots\times\mathbf{
\L}_{n_{k-1}})$. Since $\mathbf{ A}$ is critical, we have $m_{0},\ldots,m_{k-1}=1$.
Hence it satisfies condition $1)$. Suppose condition $2)$ fails, then there
exist $i\neq j$ and $s\neq r$ such that $n_{i} | n_{j}$,
$n_{s} | n_{r}$ and $j\neq r$.  Since $n_{i} | n_{j}$, we have that the
correspondence that maps
\[(a(0),\ldots,a(j-1),a(j+1),\ldots,a(k-1))\]
to
\[(a(0),\ldots,a(j-1),a(i),a(j+1),\ldots,a(k-1)).\]
defines an isomorphism from
$\mathbf{ \L}_{n_{0}}\times\cdots\times\mathbf{ \L}_{n_{j-1}}\times \mathbf{
\L}_{n_{j+1}}\times\cdots\times\mathbf{ \L}_{n_{k-1}}$ onto  a proper subalgebra of
$\mathbf{ A}$.\\
Similarly the algebra $\mathbf{ \L}_{n_{0}}\times\cdots\times \mathbf{
\L}_{n_{r-1}}\times\mathbf{ \L}_{n_{r+1}}\times\cdots\times\mathbf{ \L}_{n_{k-1}}$ is
isomorphic to a proper subalgebra of $\mathbf{ A}$. Finally, observe that $\mathbf{
\L}_{n_{0}}\times\cdots\times\mathbf{ \L}_{n_{k-1}}$ is embeddable into
\[{\mathbf{ \L}_{n_{1}}}^{2}\times\cdots\times{\mathbf{ \L}_{n_{j-1}}}^{2}\times\mathbf{
\L}_{n_{j}}\times {\mathbf{ \L}_{n_{j+1}}}^{2}\times\cdots\times{\mathbf{
\L}_{n_{r-1}}}^{2}\times\mathbf{ \L}_{n_{r}}\times{\mathbf{
\L}_{n_{r+1}}}^{2}\times\cdots\times{\mathbf{ \L}_{n_{k}}}^{2},\] by means of the
correspondence $\delta $ defined as:
 \[\delta (a(0),\ldots,a(k-1))(i) = \left\{
                                      \begin{array}{ll}
                                        (a(i),a(i)), & \hbox{if $i\not= j,r$;} \\
                                        a(i), & \hbox{if $i=j,r$.}
                                      \end{array}
                                    \right.
\]
Therefore\\
$\mathbf{ A}\in\mathbb{Q}(\mathbf{ \L}_{n_{0}}\times\cdots\times\mathbf{ \L}_{n_{j-1}}\times
\mathbf{ \L}_{n_{j+1}}\times\cdots\times\mathbf{ \L}_{n_{k-1}}$, $\mathbf{
\L}_{n_{0}}\times\cdots\times \mathbf{ \L}_{n_{r-1}}\times\mathbf{
\L}_{n_{r+1}}\times\cdots\times\mathbf{ \L}_{n_{k-1}})$ in contradiction with the
fact that $\mathbf{ A}$ is critical. \hfill $\Box$
\bg{corollary}
\label{criticalfinit} The number of non isomorphic critical MV-algebras in a
proper variety of $\mathbb{W}$ is finite.\hfill $\Box$
\e{corollary}
\textbf{Proof :} If $\mathbb{K}$ is a proper variety of $\mathbb{W}$, then it is
shown in \cite{Ko} and \cite{CDM} that  $\mathbb{K}$ contains a finite number of
${\mathbf{ \L}_{n}}'s$. Let $M=\{n\in\omega : \mathbf{ \L}_{n}\in \mathbb{K}\}$, clearly
$|M|<\omega$. By Theorem \ref{Wcritical}, all critical algebras in
$\mathbb{K}$ are:
\[\mathbb{I}(\{\mathbf{ \L}_{n_{1}}\times\cdots\times\mathbf{ \L}_{n_{l}} : \mbox{
satisfying (1) and (2) of Theorem \ref{Wcritical} and  } n_{i}\in M\}).\]
Since $|M|$ is finite we have that the number of non isomorphic
critical MV-algebras in $\mathbb{K}$ is finite.\hfill $\Box$\\
From the above result we deduce :
\bg{theorem}
\label{locfinfingen}  A quasivariety of MV-algebras is locally finite if and
only if it is finitely generated.
\e{theorem}
\textbf{Proof :} Let $\mathbb{K}$ be a locally finite quasivariety of MV-algebras,
by Corollary \ref{quasivarWlocfin}, $\mathbb{V}(\mathbb{K})$ is a proper
subvariety of $\mathbb{W}$, thus, applying Corollary \ref{criticalfinit} the
number of non isomorphic critical MV-algebras in  $\mathbb{K}$ is finite. By
Theorem \ref{Criticalgen}, $\mathbb{K}$ is generated by its critical algebras,
therefore, since a quasivariety is closed under the operation of isomorphic
images, $\mathbb{K}$ is finitely generated.\\
The converse is given by Corollary \ref{FinGenisLoc} \hfill $\Box$

\bg{corollary}
\label{finquas} Every locally finite variety of MV-algebras
contains a finite number of quasivarieties and, therefore,
the class of locally finite quasivarieties of MV-algebras is countable.\hfill $\Box$
\e{corollary}

In order to classify and distinguish locally finite quasivarieties by looking at their
generators we need the following result:

\bg{lemma}
\label{distinct}
Let $\{\mathbf{ \L}_{n_{i1}}\times\cdots\times\mathbf{ \L}_{n_{il(i)}} : i\in I\}$ and
$\{\mathbf{ \L}_{m_{j1}}\times\cdots\times\mathbf{ \L}_{m_{jl(j)}} : j\in J\}$ two finite families of
critical MV-algebras, then
\[\mathbb{Q}(\{\mathbf{ \L}_{n_{i1}}\times\cdots\times\mathbf{ \L}_{n_{il(i)}} : i\in
I\})\sub\mathbb{Q}(\{\mathbf{ \L}_{m_{j1}}\times\cdots\times\mathbf{ \L}_{m_{jl(j)}} :
j\in J\})\]
if and only if\\
For every $i\in I$, there exists  $H$, a non-empty subset of $J$, such that:
\be
 \item For any $1\le k\le l(i)$ there are $j\in H$ and $1\le r\le l(j)$ such that
   $n_{ik}|m_{jr}$.
 \item For any $j\in H$ and $1\le r\le l(j)$ there exists $1\le k\le l(i)$ such that
   $n_{ik}|m_{jr}$.
\ee
\e{lemma}
\textbf{Proof :} Assume $\mathbb{Q}(\{\mathbf{ \L}_{n_{i1}}\times\cdots\times\mathbf{
\L}_{n_{il(i)}} : i\in I\})\sub\mathbb{Q}(\{\mathbf{
\L}_{m_{j1}}\times\cdots\times\mathbf{ \L}_{m_{jl(j)}} : j\in J\})$, then for every $i\in I$, there is
$\emptyset\neq H\sub J$ such that
$\mathbf{ \L}_{n_{i1}}\times\cdots\times\mathbf{
\L}_{n_{il(i)}}$ is embeddable into $\displaystyle\prod_{j\in H} (\mathbf{
\L}_{m_{j1}}\times\cdots\times\mathbf{ \L}_{m_{jl(j)}})^{\alpha_{j}}$. Therefore, since
$\displaystyle \bigcup_{j\in H}\{m_{jr} : r\leq l(j)\}$ is a finite set
then, it follows from lemma \ref{Ln12inclusiolm12} conditions 1 and 2 .\\
To prove the converse we show that for every $i\in I$,
\[\mathbf{ \L}_{n_{i1}}\times\cdots\times\mathbf{
\L}_{n_{il(i)}}\in \mathbb{I}\mathbb{S}\mathbb{P}(\{\mathbf{ \L}_{m_{j1}}
\times\cdots\times\mathbf{ \L}_{m_{jl(j)}} : j\in H\})\]
where $H$ is the subset of $J$ defined in the hypothesis
By condition 1), for every  $1\le k\le l(i)$, we choose
$j\in H$ named $j_{k}$ and we choose  $1\le r_{k}\le l(j_{k})$
such that $n_{k}|m_{j_{k}r_{k}}$,
 the the following map
\[\beta \colon \L_{n_{i1}}\times\cdots\times
\L_{n_{il(i)}} \to \prod_{1\le k\le l(i)}\L_{m_{j_{k}1}}
\times\cdots\times\L_{m_{j_{k}l(j_{k})}}\ :\bar{a}\mapsto \beta(\bar{a})\]

\[\mbox{on }\beta(\bar{a})(k)(r)=\left\{
\bg{array}{ll}
 \bar{a}(k) & \mbox{ if } r=r_{k}\\
 \bar{a}(l)  &\bg{array}{l}
                 \mbox{for some } 1\le l\le l(i) \mbox{ such that }
                 n_{l}|m_{j_{k}r}\mbox{ if } r\neq r_{k}\\
                 \mbox{it exists by condition 2) }
               \e{array}
\e{array}\right. \]
gives an embedding. \hfill $\Box$

It is well known that every subvariety of MV-algebras is finitely axiomatizable.
In fact, some effective axiomatizations are given in \cite{D-L} and in \cite{Pa}.
To our concern, we only need to axiomatize locally finite varieties of
MV-algebras. In \cite{To} it is proved that the variety generated by
$\mathbf{ \L}_{n}$ is finitely axiomatizable and it is axiomatized  by MV1,...,MV6
plus a single axiom of the form $\varphi(x)\approx 1$, denoted by
$v_{n}(x)\approx 1$. Moreover, {\em for every $n_{1},\ldots,n_{r}<\omega , {\Bbb
V}(\mathbf{ \L}_{n_{1}},\ldots,\mathbf{ \L}_{n_{r}})$ is the subvariety of $\mathbb{W}$
defined by the equation $v_{n_{1}}(x)\lor\cdots\lor v_{n_{r}}(x)\approx 1$}
\cite[Theorem 1.8]{To}. Where  $\lor$ is defined by $x\lor
y=\lnot(\lnot x\oplus y)\oplus y$.\\

In general, locally finite quasivarieties are not finitely axiomatizable, not even finitely
generated quasivarieties are finitely axiomatizable. For instance: Let $\mathbb{K}=\mathbb{Q}(\mathbf{ A})$
where $\mathbf{ A}=\langle \{0,1,2\},f,g \rangle$ is of type (1,1) with $f$ and $g$
defined by $f(0)=1, g(0)=2$ and $f(x)=g(x)=x$ for $x\neq 0$. Due to Gorbunov
\cite{Go}, $\mathbb{K}$ is not finitely axiomati\-zable while it is finitely
generated.(see also \cite[page 149]{Cz-Dz}).In the case of MV-algebras, we will show that locally
finite quasivarieties of MV-algebras are finitely axiomati\-zable.

\bg{theorem}
\label{finitaxiomatitzable} Every locally finite  quasivariety of MV-algebras is
finitely axiomatizable.
\e{theorem}
\textbf{Proof :} Let $\mathbb{K}$ be a locally finite quasivariety of MV-algebras,
by Corollary \ref{quasivarWlocfin}, $\mathbb{V}(\mathbb{K})$ is a proper locally
finite subvariety of $\mathbb{W}$ and therefore it is finitely axiomatizable. Since $\mathbb{K}$ is
finitely generated and $\mathbb{V}(\mathbb{K})$ is finitely axiomatizable we only need to prove that
$\mathbb{K}$ can be axiomatized by a set of quasiequations with a finite number of variables
(see \cite[Lemma 2.3]{B-P00}).\\
Let $\Sigma_{n}$ be the set of quasiequations with at most $n$ variables satisfied by $\mathbb{K}$
and $Mod(\Sigma_{n})$ be the family of algebras of same type satisfying $\Sigma _{n}$. Observe
that, since every locally finite variety of MV-algebras is finitely axiomatizable with two
variables, $Mod(\Sigma _{n})$ is a quasivariety contained in $\mathbb{V}(\mathbb{K})$ for every
$n\geq 2$. Moreover, $Mod(\Sigma _{n})\supseteq Mod(\Sigma _{m})\supseteq\mathbb{K}$ for every
$n<m$. Assume that $\mathbb{K}$ is not axiomatizable with a finite number variables, then
there is an infinite sequence $n_{1}, n_{2}, n_{3}, \ldots$ such that $\mathbb{V}({\Bbb
K})\supseteq Mod(\Sigma_{n_{1}})\supsetneq Mod(\Sigma _{n_{2}})\supsetneq Mod(\Sigma
_{n_{3}})\supsetneq\ldots\supsetneq\mathbb{K}$ which contradicts the fact that $\mathbb{V}(\mathbb{K})$
contains a finite number of quasivarieties of MV-algebras (Corollary \ref{finquas}).
\hfill
$\Box$

\section{Two examples}
In order to show the difficulty and/or simplicity to obtain axiomatizations of
locally finite quasivarieties of MV-algebras we give two examples:

\bg{exam} Quasivarieties contained in ${\Bbb
V}(\mathbf{ \L}_{p}, \mathbf{ \L}_{q})$, where $p$ and $q$ are two distinct
prime natural numbers.
\e{exam}
Since the only divisors of $p$ and $q$ are $1,p,q$,
by the characterization given in Theorem \ref{Wcritical}, we have that  all
critical MV-algebras contained in
$\mathbb{V}(\mathbf{ \L}_{p}, \mathbf{ \L}_{q})$ are
\[\mathbb{I}(\{\mathbf{ \L}_{1}, \mathbf{ \L}_{p},\mathbf{ \L}_{q},\mathbf{ \L}_{1}\times\mathbf{
\L}_{p}, \mathbf{ \L}_{1}\times\mathbf{ \L}_{q}, \mathbf{ \L}_{p}\times\mathbf{ \L}_{q} \}).\]

\no Applying  Lemma \ref{distinct} we obtain all subquasivarieties of ${\Bbb
V}(\mathbf{
\L}_{p}, \mathbf{
\L}_{q})$ which are sketched in figure 1.

\bg{figure}[hbt]
\thicklines
\label{Quasivarpq}
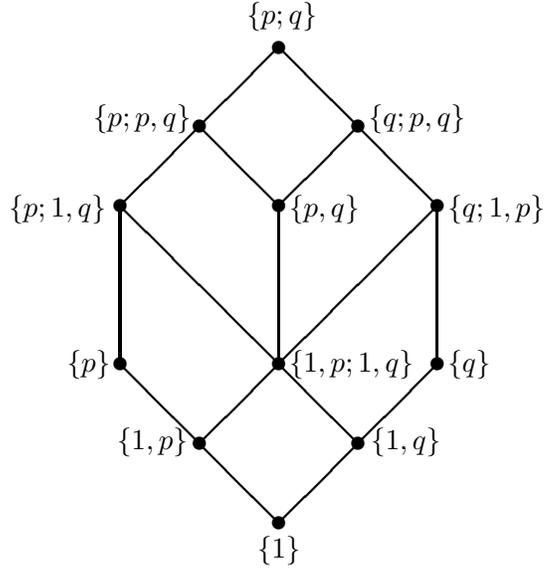
\caption{Lattice of all quasivarieties contained in $\mathbb{V}(\mathbf{
L}_{p}, \mathbf{ L}_{q})$.}
\bg{picture}(380,260)(-190,-130)
\multiput(-60,30)(60,0){3}{\line(0,-1){60}}
\multiput(0,90)(-60,-120){2}{\line(1,-1){60}}
\multiput(0,90)(60,-120){2}{\line(-1,-1){60}}
\put(-60,30){\line(1,-1){90}}
\put(60,30){\line(-1,-1){90}}
\put(0,30){\line(-1,1){30}}
\put(0,30){\line(1,1){30}}

\multiput(0,90)(0,-60){4}{\circle*{5}}
\multiput(-30,60)(0,-120){2}{\circle*{5}}
\multiput(30,60)(0,-120){2}{\circle*{5}}
\multiput(-60,30)(0,-60){2}{\circle*{5}}
\multiput(60,30)(0,-60){2}{\circle*{5}}

\put(-12,99){$\{p ; q\}$}
\put(34,60){$\{q ; p,q\}$}
\put(-70,60){$\{p ; p,q\}$}
\put(4,26){$\{p,q\}$}
\put(64,26){$\{q ; 1,p\}$}
\put(-102,26){$\{p ; 1,q\}$}
\put(64,-33){$\{q\}$}
\put(-80,-33){$\{p\}$}
\put(4,-33){$\{1,p ; 1,q\}$}
\put(35,-62){$\{1,q\}$}
\put(-61,-62){$\{1,p\}$}
\put(-8,-103){$\{1\}$}

\put(-180,-130){$\{n_{1},\ldots,n_{l};m_{1},\ldots,m_{k}\}$ stands for  ${\Bbb
Q}(\mathbf{ \L}_{n_{1}}\times\cdots\times\mathbf{ \L}_{n_{l}},\mathbf{
\L}_{m_{1}}\times\cdots\times\mathbf{ \L}_{m_{k}})$.}
\e{picture}
\e{figure}

\no Given a quasivariety $\mathbb{K}$ of MV-algebras, let
\[\mathbb{K}:\mathbf{ \L}_{n}=\{\mathbf{ A}\in\mathbb{K} : \mathbf{ \L}_{n}\not\in{\Bbb
I}\mathbb{S}(\mathbf{ A})\}.\]

We define $0x=0$ and for each $n\in\omega \ (n+1)x=x\oplus nx$. It is known
(see \cite[Lemma 2.2.]{To}) that
 $\mathbf{ \L}_{n}$ is embeddable into  an MV-algebra $\mathbf{ A}$
if, and only if, there is an element $a\in\mathbf{ A}$ such that $(n-1)(\lnot a)=a$ and
$a\neq 1^{\mathbf{ A}}$.  Thus  the quasiequation
$(n-1)(\lnot x)\approx x\Ra x\approx 1$
holds in an MV-algebra $\mathbf{A}$ if and only if $\mathbf{A}$ does not contain a copy of $\mathbf{
\L}_{n}$. Therefore, $\mathbb{K}:\mathbf{ \L}_{n}$ is the quasivariety  axiomatized by:

$$\{ \mbox{axioms of } \mathbb{K}\} \cup \{ (n-1)(\lnot x)\approx x\Ra x\approx 1 \}.$$

From the above,  with more or less difficulty, one can prove the following:

\be
 \item $\mathbb{Q}( \{\mathbf{ \L}_{p}, \mathbf{ \L}_{p}\times\mathbf{ \L}_{q}\})={\Bbb
  V}(\{\mathbf{ \L}_{p}, \mathbf{ \L}_{q}\}):\mathbf{ \L}_{q}$.
 \item $\mathbb{Q}( \{\mathbf{ \L}_{q}, \mathbf{ \L}_{p}\times\mathbf{ \L}_{q}\})={\Bbb
  V}(\{\mathbf{ \L}_{p}, \mathbf{ \L}_{q}\}):\mathbf{ \L}_{p}$.
 \item $\mathbb{Q}(\mathbf{ \L}_{p}\times\mathbf{ \L}_{q})=
 (\mathbb{V}(\{\mathbf{ \L}_{p}, \mathbf{ \L}_{q}\}):\mathbf{ \L}_{p}):\mathbf{ \L}_{q}$.
 \item $\mathbb{Q}(\mathbf{ \L}_{1}\times\mathbf{ \L}_{p})=\mathbb{V}(\mathbf{ \L}_{p}):\mathbf{
  \L}_{p}$.
 \item $\mathbb{Q}(\mathbf{ \L}_{1}\times\mathbf{ \L}_{q})=\mathbb{V}(\mathbf{ \L}_{q}):\mathbf{
  \L}_{q}$.

 \item $\mathbb{Q}(\mathbf{ \L}_{p}, \mathbf{ \L}_{1}\times\mathbf{ \L}_{q})$ is axiomatized by
  the axioms of $\mathbb{V}(\mathbf{ \L}_{p}, \mathbf{ \L}_{q}):\mathbf{ \L}_{q}$ and
  \[ rx\approx 1\ \&\  r(\lnot x)\approx 1 \Ra v_{p}(y)\approx 1,
  \mbox{ where } r=\max\{p,q\}.\]
 \item $\mathbb{Q}(\mathbf{ \L}_{q}, \mathbf{ \L}_{1}\times\mathbf{ \L}_{p})$ is axiomatized by
  the axioms of $\mathbb{V}(\mathbf{ \L}_{p}, \mathbf{ \L}_{q}):\mathbf{ \L}_{p}$ and
  \[ rx\approx 1\ \&\  r(\lnot x)\approx 1 \Ra v_{q}(y)\approx 1,\mbox{ where }
  r=\max\{p,q\}.\]
 \item $\mathbb{Q}(\mathbf{ \L}_{1}\times\mathbf{ \L}_{p}, \mathbf{ \L}_{1}\times\mathbf{
  \L}_{q})$ is axiomatized by the axioms of $(\mathbb{V}(\mathbf{ \L}_{p}, \mathbf{
  \L}_{q}):\mathbf{ \L}_{q}):\mathbf{ \L}_{p}$ and
  \[ rx\approx 1\ \&\  r(\lnot x)\approx 1 \Ra x\approx 1,\mbox{ where } r=\max\{p,q\}.\]
\ee

\bg{exam} Quasivarieties contained in $\mathbb{V}(\mathbf{ \L}_{p^{r}})$, where $p$ is a prime
natural number and $r$ is a natural number.
\e{exam}
The class of its critical algebras is given by:
\[\displaystyle \mathbb{I}\Big(\{\mathbf{ \L}_{p^{s}}:s\le r\}\cup\{\mathbf{ \L}_{p^{n}}\times\mathbf{
\L}_{p^{m}}:n<m\le r\}\Big),\]
After lengthy computations (some cases  are not straightforward),  one can prove that each
subquasivariety of
$\mathbb{V}(\mathbf{ \L}_{p^{r}})$ belongs to  one of the following three types:
\be
 \item $\mathbb{Q}(\mathbf{ \L}_{p^{s}})=\mathbb{V}(\mathbf{ \L}_{p^{s}})$ where $s\leq r$ and it is
  axiomatized by
  \be
    \item {\rm MV1},...,{\rm MV6} and
    \item $v_{p^{s}}(x)\approx 1$.
  \ee
 \item $\mathbb{Q}(\mathbf{ \L}_{p^{n_{1}}}\times\mathbf{ \L}_{p^{m_{1}}},\ldots, \mathbf{
  \L}_{p^{n_{k}}}\times\mathbf{ \L}_{p^{m_{k}}})$ such that $n_{i}<m_{i}\leq r$,
  for every $1\leq i\leq k$ and $n_{i}<n_{j}$  and $m_{i}>m_{j}$ if $i<j$ and it is
  axiomatized by:
  \be
   \item {\rm MV1},...,{\rm MV6},
   \item $v_{p^{m_{1}}}(x)\approx 1$,
   \item $(p^{n_{k}+1}-1)(\lnot x)\approx x\Ra x\approx 1$ and
   \item $(p^{n_{j-1}+1}-1)(\lnot x)\approx x\Ra v_{p^{m_{j}}}(y)\approx 1$ for
    every $2\leq j\leq k$.
  \ee
 \item $\mathbb{Q}(\mathbf{ \L}_{p^{n_{1}}}\times\mathbf{ \L}_{p^{m_{1}}},\ldots, \mathbf{
  \L}_{p^{n_{k}}}\times\mathbf{ \L}_{p^{m_{k}}}, \mathbf{ \L}_{p^{s}})$ such that
  $n_{i}<s<m_{i}\leq r$, for every $1\leq i\leq k$ and $n_{i}< n_{j}$ and
  $m_{i}>m_{j}$  if $i<j$ and it is axiomatized by:
  \be
   \item {\rm MV1},...,{\rm MV6},
   \item $v_{p^{m_{1}}}(x)\approx 1$,
   \item $(p^{s+1}-1)(\lnot x)\approx x\Ra x\approx 1,$
   \item $(p^{n_{j-1}+1}-1)(\lnot x)\approx x\Ra v_{p^{m_{j}}}(y)\approx 1$ for
    every $2\leq j\leq k$ and
   \item $(p^{n_{k}+1}-1)(\lnot x)\approx x\Ra v_{p^{s}}(y)\approx 1$.
  \ee
\ee

\no It remains an {\em open problem} whether we can find a computational method to obtain
for every locally finite quasivariety its axiomatization.

\bg{thebibliography}{99}

\bibitem{B-P00}{\sc J.Blok and D.Pigozzi},   Finitely based quasivarieties,
{\em Algebra Universalis\/} \textbf{22\/} (1986), 1-13.

\bibitem{B-P}{\sc J.Blok and D.Pigozzi},
 Algebraizable logics,   {\em Mem. Amer. Math.
Soc.\/} \textbf{396\/}, vol 77. Amer. Math. Soc.
Providence, 1989.

\bibitem{Bu}{\sc S.Burris and  H.P.Sankappanavar},
"A course in Universal Algebra,"
Springer Verlag.  New York, 1981.

\bibitem{Ch1} {\sc C.C.Chang},   Algebraic analysis of
many-valued logics, {\em Trans. Amer. Math. Soc.\/}  {\bf
88\/} (1958), 467-490.

\bibitem{Ch2} {\sc C.C.Chang},   A new proof of the
completeness of the \L ukasiewicz axioms,
{\em Trans. Amer. Math. Soc.\/} \textbf{93\/} (1959), 74-80.

\bibitem{CDM} {\sc R.Cignoli, I.M.L.D'Ottaviano  and
D.Mundici}, "Algebras das Logicas de
\L ukasiewicz," UNICAMP. Brasil, 1994.

\bibitem{C-T} {\sc R.Cignoli and  A. Torrens},
Retractive MV-Algebras, {\em Mathware \& Soft
Computing\/} \textbf{2\/} (1995), 157-165.

\bibitem{Cz-Dz}{\sc J.Czelakowski and W.Dziobiak},
Congruence distributive quasivarieties whose finitely
subdirectly irreducible members form a universal
class, {\em Algebra Universalis\/} \textbf{27\/} (1990),
128-149.

\bibitem{D-L} {\sc A.Di Nola and A.Lettieri},
Equational Characterization of all Varieties of
MV-algebra, {\em Journal of Algebra \/}\textbf{221\/}(1999), 463-474.

\bibitem{Dz0}{\sc W.Dziobiak},  On subquasivariety
lattices of semi-primal varieties, {\em Algebra
Universalis\/} \textbf{20\/} (1985), 127-129.

\bibitem{Dz1}{\sc W.Dziobiak}, Finitely generated
congruence distributive quasivarieties of algebras,
{\em Fundamenta Mathematicae\/} \textbf{133\/} (1989), 47-57.

\bibitem{Go}{\sc V.A. Gorbunov}, Covers in
subquasivariety lattices and independent
axiomatizability,  {\em Algebra
i Logika\/} \textbf{16\/} (1977), 507-548.

\bibitem{Ko} {\sc Y. Komori}, Super-\L ukasiewicz
Propositional Logic, {\em Nagoya Math. J.\/} \textbf{84\/}
(1981), 119-133.

\bibitem{Mu1} {\sc D. Mundici}, Interpretation of AF
C*-algebras in \L ukasiewicz Sentential Calculus,
{\em J. Funct. Anal.\/} \textbf{65\/} (1986),  15-63.

\bibitem{Pa} {\sc G. Panti}, Varieties of
MV-algebras, {\em Journal of Applied and Non-Classical Logic\/}
\textbf{9\/} (1999), 141-157.

\bibitem{Rod} {\sc A.J. Rodr\'iguez}, "Un Estudio
Algebraico de los C\'alculos Proposicionales de
$\L$ukasiewicz,"  Ph. D.
Thesis, University of Barcelona, 1980.

\bibitem{RTV} {\sc A.J.Rodr\'iguez, A.Torrens and
V.Verd\'u}, \L ukasiewicz logic and Wajsberg
algebras, {\em Bull. Sec. Log. Polish Ac. Sc.\/}  {\bf
2\/}, vol 19 (1990), 51-55.

\bibitem{To} {\sc A. Torrens}, Cyclic Elements in
MV-algebras and Post Algebras,  {\em Math.
Log. Quart.\/} \textbf{40\/} (1994), 431-444.

\end{thebibliography}

\noindent {\sc Joan Gispert  } and  {\sc Antoni Torrens }\\
Facultat de Matem\`atiques\\
Universitat de Barcelona\\
Gran Via de les Corts Catalanes 585\\
08007 BARCELONA (Spain)\\
{\tt $\{$gispert,torrens$\}$@mat.ub.es}

\e{document}